\documentstyle[12pt]{article}

\begin{document}
\begin  {center}{\Large {A note on Morse's index theorem for Perelman's $\mathcal{L}$-length }}
\end {center}

\begin {center} {Hong Huang}

\end {center}

 \hspace*{1cm}


 {\bf Abstract} This is essentially a note on Section 7 of Perelman's first paper
on Ricci flow. We list some basic properties of the index form for
Perelman's $ \mathcal{L} $-length, which are analogous to the ones
in Riemannian case (with fixed metric), and observe that Morse's
index theorem for Perelman's  $\mathcal{L}$-length holds. As a
corollary we get the finiteness of the number of the
$\mathcal{L}$-conjugate points along a finite
$\mathcal{L}$-geodesic.

\hspace*{0.4cm}

In his ground-breaking work [6] on Ricci flow Perelman introduced
$\mathcal{L}$-length, $\mathcal{L}$-Jacobi field among many other
important innovations. For more details see [2],[3], and [5]. Here
we'll add some notes on Section 7 of this paper of Perelman's. We
list some basic properties of the index form for Perelman's $
\mathcal{L} $-length, which are analogous to the ones in
Riemannian case (with fixed metric, cf.[1], [4] and [7] ), and
observe that Morse's index theorem for Perelman's
$\mathcal{L}$-length holds. The main idea of the proof is the same
as that of the fixed metric case, but one needs to  be careful
when the $\tau$-interval is $[0,\bar{\tau}]$ (see in particular
the proof of the Key Lemma below). As a corollary we get the
finiteness of number of the $\mathcal{L}$-conjugate points along a
finite $\mathcal{L}$-geodesic.

Throughout this note we assume $(M,g(\tau))$, where
$(g_{ij})_{\tau}=2R_{ij}$, is a (backwards) Ricci flow, which is
complete for each $\tau$-slice and has uniformly bounded curvature
operator on an interval $[\tau_{1},\tau_{2}]$.

Recall Perelman's $ \mathcal{L }$-length
$\mathcal{L}(\gamma)$=$\int_{\tau_{1}}
 ^{\tau_{2}}\sqrt{\tau}(R(\gamma(\tau))+|\dot{\gamma}(\tau)|^{2})d\tau$
for a curve $\gamma(\tau)$ ($\tau_{1}\le \tau \le \tau_{2}$) in
$M$.

 {\bf Definition 1} Let $p \in M$, $v \in T_{p}M$, and $\gamma_{v}$ be the
 $ \mathcal{L }$-geodesic with
$\gamma_{v}(0)=p $, $lim_{\tau \rightarrow
0}\sqrt{\tau}\dot{\gamma}_{v}(\tau)=v$. We say
$q=\gamma_{v}(\bar{\tau})$ is a $ \mathcal{L} $-conjugate point of
$p$ along the $\mathcal{L} $-geodesic $\gamma_{v}$ if $v$ is a
critical point of the $\mathcal{L}$$_{\bar{\tau}}$-exponential map
$\mathcal{L}$$ _{\bar{\tau}}$$exp$. Here, (following the notation
in [2],)
$\mathcal{L}$$_{\bar{\tau}}$$exp(v)$=$\mathcal{L}$$exp_{v}(\bar{\tau})$
(=$\gamma_{v}(\bar{\tau})$).

{\bf Definition 2}(Perelman) A vector field along a  $\mathcal{L}
$-geodesic $\gamma$ is called
 $\mathcal{L}$-Jacobi field, if it is the variation vector field of a one parameter
 family of $\mathcal{L}$-geodesics $\gamma_{s}$ with $\gamma_{0}=\gamma$ .

The equation for a $\mathcal{L}$-Jacobi field $U$ along a
$\mathcal{L} $-geodesic $\gamma_{v}(\tau)$ ( $0 < \tau_{1}\le \tau
\le \tau_{2}$) is ( see, for example, [2])

$\nabla_{X}\nabla_{X}U-R(X,U)X-1/2\nabla_{U}(\nabla
R)+2(\nabla_{U}$Ric)$(X)+2$Ric$(\nabla_{X}U)+1/(2\tau)\nabla_{X}U=0.$

(Here and below, $X(\tau)=\dot{\gamma}_{v}(\tau)$. Moreover
$Ric(Y)$ here means $Ric(Y,\cdot)$ in Perelman [6].   )

One can easily extend this to the case $\tau_{1}=0$. (See [2].)

{\bf Remark 1} As in the Riemannian case (with fixed metric)
$q=\gamma_{v}(\bar{\tau})$ is a $ \mathcal{L} $-conjugate point of
$p=\gamma_{v}(0)$ along the $\mathcal{L} $-geodesic
$\gamma_{v}(\tau)$($0 \le \tau \le \bar{\tau}$) if and only if
there is a nontrivial $\mathcal{L}$-Jacobi field $U$ along
$\gamma_{v}$ with $U(0)=U(\bar{\tau})=0$.

The $\mathcal{L}$-index form along a $\mathcal{L}$-geodesic
$\gamma_{v}(\tau)$ ( $ \tau_{1}\le \tau \le \tau_{2}$) is defined
to be

$I(U,V)=\int_{\tau_{1}} ^{\tau_{2}}\tau^{1/2}$[(Hess$
R)(U,V)+2\langle \nabla_{X}U,\nabla_{X}V \rangle+2\langle
R(U,X)V,X\rangle-2(\nabla_{U}$Ric)$(V,X)-2(\nabla_{V}$Ric)$(U,X)+2(\nabla_{X}$Ric)$(U,V)]d\tau$

for any piecewise smooth vector fields $U$,$V$ along $\gamma_{v}$.

This is a symmetric, bilinear form.

 {\bf Remark 2} If $Y(\tau)$ is a smooth
vector field along a $\mathcal{L}$-geodesic $\gamma_v(\tau) (0 \le
\tau \le \bar{\tau})$, and $Y(0)=0$, then
$I(Y,Y)=\delta_{Y}^{2}$$\mathcal{L}$$-\delta_{\nabla_{Y}Y}$$\mathcal{L}$.
(Compare with formula (7.7) in Perelman [6].)

Now we prove a

{\bf Key Lemma} For any  vector fields $U$,$V$ along a
$\mathcal{L}$-geodesic $\gamma_{v}(\tau)$($0< \tau_{1}\le \tau \le
\tau_{2}$) with $U$ smooth ( and $V$ piecewise smooth ), we have

$I(U,V)=2\tau^{1/2}\langle \nabla_{X}U,V\rangle|_{\tau_{1}}
^{\tau_{2}}-2\int_{\tau_{1}} ^{\tau_{2}}\tau^{1/2}\langle
\nabla_{X}\nabla_{X}U-R(X,U)X-1/2\nabla_{U}(\nabla
R)+2(\nabla_{U}$Ric)$(X)+2$Ric$(\nabla_{X}U)+1/(2\tau)\nabla_{X}U,V\rangle
d\tau.$

Furthermore, this equality extends to the case $\tau_{1}=0$.

{\bf Proof} In case $\tau_{1}>0$ one simply use

 $d/d\tau \langle
\nabla_{X}U,V\rangle=\langle
\nabla_{X}\nabla_{X}U,V\rangle+\langle \nabla_{X}U,\nabla_{X}V
\rangle+2$Ric$(\nabla_{X}U,V)+(\nabla_{X}$Ric)$(U,V)+(\nabla_{U}$Ric)$(V,X)-(\nabla_{V}$Ric)$(X,U)$

(compare with formula (11.2) in [3]), and

$d/d\tau (\tau^{1/2}\langle
\nabla_{X}U,V\rangle)=(1/2)\tau^{-1/2}\langle
\nabla_{X}U,V\rangle+\tau^{1/2}d/d\tau \langle
\nabla_{X}U,V\rangle$,

 then integration by parts, and the formula follows.

 To justify the $\tau_{1}=0$ case, it suffices to  observe

$\tau(\nabla_{X}\nabla_{X}U-R(X,U)X-1/2\nabla_{U}(\nabla
R)+2(\nabla_{U}$Ric)$(X)+2$Ric$(\nabla_{X}U)+1/(2\tau)\nabla_{X}U)
=\nabla_{\sqrt{\tau}X}\nabla_{\sqrt{\tau}X}U-R(\sqrt{\tau}X,U)\sqrt{\tau}X
-\tau/2\nabla_{U}(\nabla
R)+2\sqrt{\tau}(\nabla_{U}$Ric)$(\sqrt{\tau}X)+2\sqrt{\tau}$Ric$(\nabla_{\sqrt{\tau}X}U)$

(compare with [2]), and note that $lim_{\tau \rightarrow
0}\sqrt{\tau}X(\tau)$ exists, and that the integration
$\int_{0}^{1}\tau^{-1/2}d\tau$ converges.

\hspace*{0.4cm}

 {\bf Remark 3} One can easily generalize the Key Lemma to the
 case that $U$ is piecewise smooth.

 \hspace*{0.4cm}

Below we list some basic properties of $\mathcal{L}$-Jacobi field
which is analogous to the Riemannian case (with fixed metric, see
[1], [4], and in particular [7]), whose proof is similar to the
fixed metric case and is omitted ( in the proof the Key Lemma
above play an important role).

For convenience we denote by $\mathcal{V}$$_0 (\tau_{1},\tau_{2})$
the space of piecewise smooth vector fields $V(\tau)$ along a
$\mathcal{L}$-geodesic $\gamma_{v}(\tau)$ ( $\tau_{1}\le \tau \le
\tau_{2})$ with $V(\tau_{1})=V(\tau_{2})=0.$

In the following lemmata we suppose $\gamma_{v}(\tau)$ (
$\tau_{1}\le \tau \le \tau_{2})$ is a $\mathcal{L}$-geodesic.

 \hspace*{0.4cm}

 {\bf  Lemma 1} Let $\gamma_{v}(\tau_{2})$ be $\mathcal{L}$-conjugate to $
\gamma_{v}(\tau_{1})$ along $\gamma_{v}$. Then for any
$\mathcal{L}$-Jacobi field $U \in $ $\mathcal{V}$$_0
(\tau_{1},\tau_{2})$ one has $I(U,U)=0$.

\hspace*{0.4cm}

 {\bf Lemma 2} If $\gamma_v(\tau)$(
$\tau_{1}\le \tau \le \tau_{2}) $ does not contain any
$\mathcal{L}$-conjugate point of $\gamma_v(\tau_{1})$, then the
$\mathcal{L}$-index form is positive definite on $\mathcal{V}$$_0
(\tau_{1},\tau_{2})$.

\hspace*{0.4cm}

{\bf  Lemma 3} Let $\gamma_{v}(\tau_{2})$ be
$\mathcal{L}$-conjugate to $ \gamma_{v}(\tau_{1})$ along
$\gamma_{v}$, but for any $\tau$ such that $\tau_{1}< \tau<
\tau_{2}$, $\gamma_{v}(\tau)$ is not $\mathcal{L}$-conjugate to $
\gamma_{v}(\tau_{1})$ along $\gamma_{v}$. Then the
$\mathcal{L}$-index form is positive semi-definite ( but not
positive definite) on $\mathcal{V}$$_0 (\tau_{1},\tau_{2})$.

\hspace*{0.4cm}

 {\bf  Lemma 4} There exists $\tau'$ with $\tau_{1}< \tau' < \tau_{2}$
 such that $\gamma_{v}(\tau')$ is $\mathcal{L}$-conjugate to $
\gamma_{v}(\tau_{1})$ along  $\gamma_{v}$ if and only if there
exists a vector field $Y\in
 \mathcal{V}$$_0 (\tau_{1},\tau_{2})$ such that $I(Y,Y)<0$.

\hspace*{0.4cm}

{\bf  Lemma 5} $U$ is a $\mathcal{L}$-Jacobi field if and only if
$I(U,Y)=0$ for any vector field $Y\in
 \mathcal{V}$$_0 (\tau_{1},\tau_{2})$ .

\hspace*{0.4cm}

{\bf  Lemma 6} Suppose $\gamma_v(\tau)$ does not contain any
$\mathcal{L}$-conjugate point of $\gamma_v(\tau_{1})$ . Let $U,Y $
be piecewise smooth vector field along $\gamma_v(\tau)$ with
$U(\tau_{1})=Y(\tau_{1})$, $U(\tau_{2})=Y(\tau_{2})$, and $U$ is a
$\mathcal{L}$-Jacobi field. Then $I(U,U)\le I(Y,Y)$. The equality
holds if and only if $Y=U$.

\hspace*{0.4cm}

 {\bf Remark 4} Lemma 6 was used in Perelman [6](7.11).

\hspace*{0.4cm}

 {\bf Lemma 7} Suppose that
$\gamma_v(\tau_{2})$ is not $\mathcal{L}$-conjugate to
$\gamma_v(\tau_{1})$ along $\gamma_v$. Then given any $w\in
T_{\gamma_v(\tau_{2})}M$ there exists an unique
$\mathcal{L}$-Jacobi  field $U$ along $\gamma_v$ such that
$U(\tau_{1})=0$ and $U(\tau_{2})=w$.

 \hspace*{0.4cm}

Let $\gamma_{v}(\tau)(0\le \tau \le \bar{\tau})$ be a
$\mathcal{L}$-geodesic. The index of the $\mathcal{L}$-index form
$I$ along $\gamma_{v}$ is defined to be the maximum dimension of a
subspace of $\mathcal{V}$$_0 (0,\bar{\tau})$ on which $I$ is
negative definite.

Now we can state Morse's index theorem for Perelman's
$\mathcal{L}$-length.

 {\bf Theorem} The index of $\mathcal{L}$-index form  along a
  $\mathcal{L}$-geodesic $\gamma_v(\tau) (0 \le \tau \le
\bar{\tau})$ is equal to the number (counting with multiplicity)
of
 $\mathcal{L}$-conjugate points $\gamma_v(\tau')$($0 < \tau' < \bar{\tau}$)
  of $\gamma_{v}(0)$ along $\gamma_v $. The index is always
  finite.

Here, by definition, the multiplicity of a $\mathcal{L}$-conjugate
point $\gamma_v(\tau')$ of $\gamma_v(0)$ along a
$\mathcal{L}$-geodesic $\gamma_v$ is the dimension of subspace
that consists of all $\mathcal{L}$-Jacobi fields in
$\mathcal{V}$$_0 (0,\tau')$.

 {\bf Proof} As in the fixed metric case, the main idea is trying to
 reduce the problem to a finite dimensional subspace of $\mathcal{V}$$_0
 (0,\bar{\tau})$, using the lemmata above. The detail is similar to that of the fixed metric case (
 see [1], [4] and [7]) and is omitted.

\hspace*{0.1cm}

 We have the following immediate

 {\bf Corollary } The number of the
$\mathcal{L}$-conjugate points of $\gamma_{v}(0)$ along
$\gamma_{v}(\tau) (0 \le \tau \le \bar{\tau})$ is finite.

\hspace*{0.1cm}

\hspace*{0.1cm}

{\bf Acknowledgements} I'm partially supported by a fund from
Beijing Normal University. I would also like to thank Prof.
Hongzhu Gao for his support.

\hspace*{0.4cm}

{\bf References}

[1]J. Cheeger, D. Ebin, Comparison theorems in Riemannian
geometry, North-Holland Publishing Co. (1975).

[2]B. Chow, P. Lu, L. Ni, A quick introduction to Ricci flow, book
to appear.

[3]B. Kleiner, J. Lott, Notes on Perelman's papers, December
30,2004.

[4]J. Milnor, Morse theory, Princeton University Press (1963).

[5]N. Sesum, G. Tian, and X. Wang, Note on Perelman's paper on the
entropy formula for the Ricci flow and its geometric applications,
October 7, 2004.

 [6]G. Perelman, The entropy formula for the Ricci
flow and its geometric applications, axXiv:math.DG/0211159 v1 11
Nov 2002.

[7]H.Wu, C. Shen and Y. Yu, Introduction to Riemannian geometry (
in Chinese), Peking University Press, 1989.

 \vspace*{1cm}

Department of Mathematics, Beijing Normal University, Beijing
100875,  P.R. China

huanghong74@yahoo.com.cn

 \end{document}